\documentclass{article}
\usepackage{amsmath,amssymb,amsthm,amscd, mathtools, latexsym}
\usepackage{enumerate,varioref, dsfont, fancyhdr}
\usepackage{tikz-cd}
\usepackage{enumitem}
\usepackage{multicol}
\usepackage{hyperref}
\usepackage{url}
\usepackage{textcomp}

\newtheorem{Thm}{Theorem}[section]
\newtheorem{Lem}[Thm]{Lemma}
\newtheorem*{Lem*}{Lemma}

\newtheorem{Not}[Thm]{Notation}

\newtheorem{Cons}[Thm]{Construction}
\newtheorem{Rem}[Thm]{Remark}

\newtheorem{Def}[Thm]{Definition}

\newtheorem{Ass}[Thm]{Assumption}
\newtheorem{Exam}[Thm]{Example}

\newcommand\restr[2]{\ensuremath{\left.#1\right|_{#2}}}

\def\cit{{\mathbb C}}

\def\qit{{\mathbb Q}}
\def\zit{{\mathbb Z}}

\def\git{{\mathbb G}}
\def\pit{{\mathbb P}}

\def\0{{\mathcal O}}

\def\C{{\mathcal C}}

\def\E{{\mathcal E}}
\def\W{{\mathcal W}}

\usepackage[title]{appendix}

\begin{document}
\title{On the integral Hodge conjecture for varieties with trivial Chow group}
\author{Humberto A. Diaz}
\newcommand{\Addresses}{{\bigskip \footnotesize
\textsc{Department of Mathematics, Washington University, St. Louis, MO 63130} \par \nopagebreak
\textit{Email address}: \ \texttt{humberto@wustl.edu}}}

\date{}
\maketitle

\begin{abstract}
\noindent We obtain examples of smooth projective varieties over $\cit$ that violate the integral Hodge conjecture and for which the total Chow group is of finite rank. Moreover,  we show that there exist such examples defined over number fields.
 
\end{abstract}

\section{Introduction}
\noindent In this note,  we consider smooth projective varieties over $\cit$ violating the integral Hodge conjecture; i.e.,  varieties $X$ with some $n$ for which there exists some $\alpha \in H^{2n} (X,  \zit(n))$ that is not algebraic but for which some multiple of $\alpha$ is.  Violations of the integral Hodge conjecture are known to exist in all degrees $2n$ for $n\neq 0,  1,  dim(X)$ and have been known to exist since work of Atiyah-Hirzebruch \cite{AH}.  Violations satisfying additional hypotheses have been noted by a number of different authors, including some rather recent ones (see,  for instance,  \cite{Ko},  \cite{SV},  \cite{CTV},  \cite{T}, \cite{BO}, \cite{CT},  \cite{OS},  \cite{Schr}).  \\
\indent In particular,  there is the recent result of Ottem-Suzuki \cite{OS},  which shows that there are examples of smooth projective varieties of dimension $3$ violating the Hodge conjecture in degree $4$ for which $CH_{0} (X) \cong \zit$ (a similar example for which the latter assumption holds conjecturally can also be found in \cite{CTV} \S 5.7).  Motivated by this,  one can consider whether there exist violations of the integral Hodge conjecture involving smooth projective varieties,  all of whose Chow groups are as small as possible.  More precisely,  we consider whether there are violations of the integral Hodge conjecture for which the total Chow group $CH^{*} (X)_{\qit}$ is of finite $\qit$-rank.  We refer to such varieties as {\em Chow-trivial} varieties.  (It is a consequence of Roitmann's Theorem \cite{B} that any such $X$ will satisfy $CH_{0} (X) \cong \zit$.) To this end,  we have the following result:
\begin{Thm}\label{main} For all $d\geq 4$,  there exist smooth projective varieties $V$ over $\cit$ of dimension $d$ that are Chow-trivial and which violate the integral Hodge conjecture in degree 4; more precisely,  for which there exists a non-algebraic torsion class in $H^{4} (V,  \zit(2))$.  In fact, one can find $V$ of this type that are defined over a number field.
\end{Thm}
\noindent Several remarks should be made at the outset.  First,  it is sufficient to consider the case that $V$ is of dimension $4$.  Once the result is proved in this case,  one can obtain results in all higher dimensions $d$ in the standard way; i.e., by replacing $V$ with $V \times \pit^{d-4}$ and using the projective bundle formula. \\ 
\indent Additionally,  our example in dimension $4$ has the form $S_{1} \times S_{2}$,  where $S_{1}$ and $S_{2}$ are Enriques surfaces.  The method of proof we use is the same method of proof employed by Colliot-Th\'el\`ene in \cite{CT},  which combines the main result of \cite{CTV} with a degeneration method of Gabber \cite{G}.  This kind of argument works best on examples involving elliptic curves,  and so we show in the sequel how the problem can be reduced to one involving elliptic curves.  In our approach,  the degeneration part of this argument involves reduction to positive characteristic,  and we explain in the final section why this is necessary.  We also observe that Schreieder has recently generalized the method of \cite{CT} to obtain violations in higher degree using refined unramified cohomology  \cite{Schr}. \\
\indent Finally,  our note leaves open the question of whether there are violations to the integral Hodge conjecture involving {\em three}-dimensional Chow-trivial varieties.  The primary difficulty is that verifying Chow-triviality in cases where it is expected is difficult in general.
\section*{Acknowledgments}
The author would like to thank Chuck Doran for useful conversations and for providing several references,  as well as Matt Kerr for his interest.

\section{Chow-trivial varieties}\label{Bloch}
We begin with the following well-known result relating the various notions of triviality of the Chow group:
\begin{Thm}\label{equiv} Let $X$ be a smooth projective variety over $\cit$.  Then,  the following are equivalent:
\begin{enumerate}[label=(\alph*)]
\item\label{triv} $X$ is Chow-trivial;
\item\label{univ-triv} for all algebraically closed extensions $K/\cit$,  $CH^{*} (X \times_{\cit} K)_{\qit}$ is of finite $\qit$-rank;
\item\label{num} rational equivalence and numerical equivalence ($\otimes \qit$) coincide for algebraic cycles on $X$;
\item\label{mot} the rational Chow motive $M(X)_{\qit}$ is a sum of Lefschetz motives;
\item\label{cyc} the total cycle class map to singular cohomology:
\[ CH^{*}(X^{m})_{\qit} \to H^{2*} (X^{m}, \qit(*))\] 
is an isomorphism for all $m$-fold powers $X^{m} = X \times \ldots \times X$.
\end{enumerate}
\begin{proof} This follows from Theorem 4 of \cite{Via}.
\end{proof}
\end{Thm}

\noindent Note that \ref{cyc} in Theorem \ref{equiv} shows that Chow-trivial varieties automatically satisfy the usual Hodge conjecture.  So,  it is natural to wonder if these also satisfy the integral Hodge conjecture.  As a result,  we obtain the following examples of Chow-trivial varieties over $\cit$:
\begin{Exam}\label{ex}
\begin{enumerate}[label=(\alph*)]
\item(Cellular varieties) If $X$ admits a stratification with strata isomorphic to a disjoint union of $\mathbb{A}^{n}$'s, then $X$ is Chow-trivial.  This follows from an argument involving the localization sequence for Chow groups.
\item(Products) If $X$ and $Y$ are Chow-trivial, then so is $X \times Y$.  This follows from characterization \ref{mot} in Theorem \ref{equiv}.
\item(Projective bundles) If $X$ is Chow-trivial, then so is any projective bundle over $X$.  This follows from a projective bundle formula argument.
\item(Enriques surfaces) Any Enriques surface is Chow-trivial (by the main result of \cite{BKL}).
\end{enumerate}
\end{Exam}
\begin{Rem}
Since our examples involve surfaces,  we note that it is classical that if $X$ is a smooth projective surface over $\cit$ that is Chow trivial, then $h^{1} (X) = 0 = p_{g}(X)$ \cite{Mum}.  Bloch's conjecture predicts that the converse also holds; i.e.,  if $X$ is a smooth projective surface over $\cit$ with $h^{1} (X) = 0 = p_{g}(X)$,  then, $X$ is Chow trivial.  Bloch's conjecture is known for all surfaces not of general type and has been verified for a number of surfaces of general type; for instance,  see \cite{GP},  \cite{V} and \cite{PW}.  It is possible that the method we use to prove Theorem \ref{main} would also work on examples involving these surfaces.
\end{Rem}
\section{Unramified cohomology}
\begin{Def} Let $X$ be a smooth variety over $\cit$, $G$ be a finite Abelian group and $\mathcal{H}_{X}^{i}(G)$ be the Zariski sheaf over $X$ associated to the presheaf $U \mapsto H^{i} (U, G)$,  where $H^{*}$ denotes singular cohomology. Then,
\[ H^{i}_{nr} (X, G) := H^{0} (X, \mathcal{H}_{X}^{i}(G))  \]
is the {\em $i^{th}$ unramified cohomology group with coefficients in $G$}.
\end{Def}
\noindent Unramified cohomology was defined in this way in \cite{BlO} (see also \cite{CTO}).  The following facts are all well-known (and either in or derivable from these two sources):
\begin{enumerate}[label=(\alph*)]
\item\label{res} $H^{i}_{nr} (X,  G)$ is a subgroup of $H^{i}(\cit(X),  G)$,  where the latter group denotes the Galois cohomology of $\cit(X)$.
\item\label{com} There is a natural map $H^{i} (X,  G) \to H^{i}(\cit(X),  G)$ factoring  as $H^{i} (X,  G) \to H^{i}_{nr} (X,  G) \hookrightarrow H^{i}(\cit(X),  G)$ whose kernel is $N^{1}H^{i} (X,  G)$,  where $N^{*}$ denotes the coniveau filtration on cohomology.
\item\label{trivial} $H^{1} (X,  G) \to H^{1}_{nr} (X,  G)$ is an isomorphism and $H^{2} (X,  G) \to H^{2}_{nr} (X,  G)$ is surjective; for $i>2$,  not much is known about $H^{i} (X,  G) \to H^{i}_{nr} (X,  G)$ apart from some special cases.
\end{enumerate}
\noindent The work of Colliot-Th\'el\`ene and Voisin \cite{CTV} proved the relation between unramified cohomology in degree $3$ and the integral Hodge conjecture in degree $4$.  This work has since been generalized by Schreieder \cite{Schr} to a relation between refined unramified cohomology and the integral Hodge conjecture in higher degree.
\begin{Thm}\label{key} When $X$ is Chow-trivial,  the integral Hodge conjecture for $X$ in degree $4$ is false $\Leftrightarrow$ $H^{3}_{nr} (X,  \zit/m(2)) \neq 0$ for some some $m$.  Moreover,  if $\gamma \in H^{3} (X,  \zit/m(2))$ has nonzero image in $H^{3} (\cit(X),  \zit/m(2))$,  then $\delta(\gamma) \in H^{4} (X,  \zit(2))$ is a non-algebraic $m$-torsion class,  where $\delta: H^{3} (X,  \zit/m(2))\to H^{4} (X,  \zit(2))$ is the coboundary in the long exact sequence in cohomology associated to the short exact sequence $0 \to \zit \to \zit \to \zit/m \to 0$.
\begin{proof} The first statement follows directly from \cite{CTV} Th\'eor\`eme 3.9.  For the second statement,  let $Z^{4} (X)$ be the quotient of $H^{4} (X,  \zit(2))$ by the subgroup of algebraic classes.  Then,  the second statement follows directly from the commutativity of the diagram:
\[\begin{tikzcd}
H^{3} (X,  \zit/m(2)) \arrow{r}{\delta} \arrow{d} & H^{4} (X,  \zit(2))[m] \arrow{d}\\
H^{3}_{nr} (X,  \zit/m(2)) \arrow{r}{\cong} & Z^{4} (X)[m]
\end{tikzcd}
\]
The commutativity follows from the fact that the isomorphism in the bottom row arises from the coboundary map in the analogous sequence for unramified cohomology.
\end{proof}
\end{Thm}
\section{Proof of Theorem \ref{main}}
As noted earlier,  it will suffice to prove that we can find such a $V$ of dimension $d=4$.  Moreover,  we will take $V$ to be a product of two Enriques surfaces defined over number fields.  Section \ref{Bloch} shows that $V$ is Chow-trivial.  Now, we let $S$ denote an Enriques surface satisfying the following criterion:
\begin{Ass}\label{ass} $S = Y/\phi$,  where $Y$ is a $K3$ surface and $\phi$ is an Enriques(=fixed-point-free) involution on $Y$,  both defined over a number field $K$.  Moreover,  $\phi$ acts on a genus $1$ curve $\tilde{E} \subset Y$,  also defined over $K$ and with $\tilde{E}(K) \neq \emptyset$; let $\pi: \tilde{E} \to E:=\tilde{E}/\phi$ be the induced quotient.  Assume that there is some place $\mathfrak{p}$ of $K$, not lying over $2$,  over which $\tilde{E}$ has multiplicative reduction.
\end{Ass}
\noindent Moreover, let $S_{2}$ denote an Enriques surface defined over $K$ with good reduction at a prime $\mathfrak{p}$ lying over $p$ (same $\mathfrak{p}$ as in Assumption \ref{ass}).  Note that for a particular $p$,  it may be nontrivial to show that such an Enriques surface exists,  so we do not claim that what follows works for any odd prime $p$ but only that it works for any odd prime for which there exists such an Enriques surface.  (Indeed,  there must be some odd primes for which this is the case, since there are certainly Enriques surfaces defined over $K$ which will necessarily have good reduction at almost all primes.) Then,  take $V = S \times S_{2}$ and let
\[ \gamma:= \pi_{S}^{*}\alpha\cup\pi_{S_{2}}^{*}\beta \in  H^{3} (V, \zit/2(2))\]
where $\alpha \in H^{1} (S, \zit/2(1)) \cong H^{1} (S, \zit/2)$ is the class corresponding to the double cover $Y \to S$ and $\beta \in H^{2} (S_{2}, \zit/2(1))$ is any class mapping nontrivially via the map:
\begin{equation}H^{2} (S_{2}, \zit/2(1)) \cong H^{2}_{\text{\'et}} (S_{2,  \cit}, \zit/2(1)) \to Br(S_{2,  \cit})[2] \cong \zit/2 \label{comp} \end{equation}
where the first arrow is the canonical isomorphism between singular and \'etale cohomology and the second arrow is from the Kummer sequence.  By Theorem \ref{key},  it will suffice to show that the image of $\gamma$ under $H^{3} (V, \zit/2(2)) \to H^{3} (\cit(V), \zit/2(2))$ is non-zero or, equivalently, $\gamma \not\in N^{1}H^{3} (V, \zit/2(2))$.  Now,  by the functoriality of the coniveau filtration (there are many ways to see this; e.g.,  see \cite{PR} \S 2.1) it will suffice to show that
$\gamma' = \restr{\gamma}{W} \not\in N^{1}H^{3} (W,  \zit/2(2))$
where $W=E \times S_{2}$ and $E \subset S$ is as above.  For the remainder of the proof,  we use \'etale cohomology with $\zit/2$ coefficients.  In particular,  we may view: 
\begin{equation} \begin{split}
\alpha':=\restr{\alpha}{E} \in H^{1}_{\text{\'et}} (E_{\cit}, \zit/2(1)) \cong H^{1}_{\text{\'et}} (E_{\overline{K_{\mathfrak{p}}}}, \zit/2(1)),  \\ 
\beta  \in   H^{2}_{\text{\'et}} (S_{2,  \cit}, \zit/2(1)) \cong H^{2}_{\text{\'et}} (S_{2,  \overline{K_{\mathfrak{p}}}}, \zit/2(1))\end{split}\label{alph} \end{equation}
where we view $E/K_{\mathfrak{p}}$ as a genus $1$ curve over $K_{\mathfrak{p}}$ (and similarly view $S/K_{\mathfrak{p}}$ as an Enriques surface over $K_{\mathfrak{p}}$) and use the standard base change results from \'etale cohomology (see \cite{M} \S VI.2).  Let $W = E \times S_{2}/K_{\mathfrak{p}}$,  again viewed as a variety over $K_{\mathfrak{p}}$.  We need to show that $\gamma' \in H^{3}_{\text{\'et}} ( W_{\overline{K_{\mathfrak{p}}}}, \zit/2(2))$ does not vanish under the map to the Galois cohomology:
\begin{equation} H^{3}_{\text{\'et}} (W_{\overline{K_{\mathfrak{p}}}}, \zit/2(2)) \to H^{3} (\overline{K_{\mathfrak{p}}}(W), \zit/2(2))\label{func} \end{equation}
\noindent To this end,  we will need to spread out the classes in (\ref{alph}).  In particular,  we will need the following lemma:
\begin{Lem}\label{Neron} Let $\mathfrak{o}$ denote the corresponding ring of integers of $K_{\mathfrak{p}}$.  There are smooth schemes $\E, \tilde{\E} \to \text{Spec} (\mathfrak{o})$ whose generic fibers are $E$ and $\tilde{E}$ and whose special fibers, $E_{0}$ and $\tilde{E}_{0}$, are $k$-isomorphic to $\git_{m, k}$.  Moreover,  $\pi: \tilde{E} \to E$ extends to an \'etale morphism $\pi: \tilde{\E} \to \E$ over $\text{Spec} (\mathfrak{o})$. 
\begin{proof}[Proof of Lemma] $\tilde{E}(K)\neq \emptyset$ by assumption,  so the Neron model of $\tilde{E}$,  $\tilde{\E} \to \text{Spec} (\mathfrak{o})$, has the structure of an Abelian group scheme. Since $\phi$ acts on $ \tilde{E}$ by translation by a two-torsion point (since $\phi$ acts freely on $\tilde{E}$),  this action extends to an action on $\tilde{\E} \to \text{Spec} (\mathfrak{o})$. We denote the induced quotient by $\E \to \text{Spec} (\mathfrak{o})$ and note that there is a corresponding \'etale double cover $\tilde{\E} \to \E$.  Finally,  we observe that since $\tilde{E}$ has multiplicative reduction over $\mathfrak{p}$,  so does the special fiber $E_{0}$ of $\E \to \text{Spec} (\mathfrak{o})$.  Hence,  the second statement.
\end{proof}
\end{Lem}
\indent Now,  we let $\mathcal{E} \to \text{Spec} (\mathfrak{o})$ be the corresponding model from Lemma \ref{Neron} and let $\mathcal{S}_{2} \to \text{Spec} (\mathfrak{o})$ denote a smooth projective morphism whose generic fiber is $S_{2}/K_{\mathfrak{p}}$.  Then,  consider $\W:= \E \times_{\mathfrak{o}}\mathcal{S}_{2}$,  whose generic fiber is $W$.
\begin{Lem}\label{lift} There is some finite unramified Galois extension $L/K_{\mathfrak{p}}$ with ring of integers $\mathfrak{o}'$ and residue field $k'$ for which there exist $A' \in H^{1}_{\text{\'et}} (\E \times_{\mathfrak{o}}\mathfrak{o}' , \zit/2(1))$ and $B \in H^{2}_{\text{\'et}} (\mathcal{S}_{2}\times_{\mathfrak{o}}\mathfrak{o}', \zit/2(1))$ lifting $\alpha'$ and $\beta$,  resp.
\begin{proof} We can describe the first class explicitly.  First note that one can take $A' \in H^{1}_{\text{\'et}} (\E, \zit/2(1))$ to be the class of the \'etale double cover (over $\mathfrak{o}$) $\tilde{\E} \to \E$ in $H^{1}_{\text{\'et}} (\E, \zit/2(1)) \cong H^{1}_{\text{\'et}} (\E, \zit/2)$.  Indeed, since $2$ is invertible in the residue field $k$ of $\mathfrak{o}$, this group parametrizes $\zit/2$-torsors over $\E$ (i.e., \'etale double covers of $\E$) (see,  for instance,  \cite{M} Chapter III).  Moreover,  the class of $\alpha'$ in (\ref{alph}) is nothing but the class of the \'etale double cover $\tilde{E}_{\overline{K_{\mathfrak{p}}}} \to E_{\overline{K_{\mathfrak{p}}}}$,  which is the fiber product of $\tilde{\E} \to \E$ with a geometric point $\text{Spec} (\overline{K_{\mathfrak{p}}}) \to \text{Spec} (\mathfrak{o})$.  This then gives a lift of $\alpha'$.  Note that for this class,  we did not need to extend $K_{\mathfrak{p}}$ yet.\\
\indent To lift $\beta$,  we set $\mathcal{S}':= \mathcal{S}_{2}\times_{\mathfrak{o}}\mathfrak{o}'$ and $S'$ the generic fiber. We need to find a suitable finite unramified Galois extension $L/K_{\mathfrak{p}}$ with ring of integers $\mathfrak{o}'$ for which the base extension map:
\[ Br(\mathcal{S}')[2] \to Br(S'_{\overline{K_{\mathfrak{p}}}})[2] \cong \zit/2\]
is surjective.  To this end,  note that the special fiber $S_{0}'$ of $\mathcal{S}' \to \text{Spec} (\mathfrak{o}')$ is smooth and projective by assumption.  So,  we observe that there is a commutative diagram:
\begin{equation} \begin{tikzcd}
Br(\mathcal{S}')[2] \arrow{r} \arrow[two heads]{d}& Br(S'_{\overline{K_{\mathfrak{p}}}})[2] \arrow{d}{\cong}\\
Br(S_{0}')[2] \arrow{r}  & Br( S_{0, \overline{k}}')[2] 
\end{tikzcd}\label{CD}\end{equation}
where the horizontal arrows are the base extension maps,  the left vertical arrow is the restriction map to the special fiber, which is surjective (by \cite{M} Chapter VI Cor.  2.7 and the Kummer sequence),  and the right vertical arrow is induced by the specialization isomorphism (\cite{M} Chapter VI Cor.  4.2).  Indeed,  there is a specialization isomorphism:
\[ H^{2}_{\text{\'et}} (S'_{\overline{K_{\mathfrak{p}}}},  \zit/2(1)) \xrightarrow{\cong} H^{2}_{\text{\'et}} (S_{0, \overline{k}}',  \zit/2(1)) \]
which maps $Pic(S_{\overline{K_{\mathfrak{p}}}})/2$ injectively into $Pic(S_{0, \overline{k}}')/2$.  Since both these groups have the same $\zit/2$-rank,  it follows that the induced map \[Pic(S'_{\overline{K_{\mathfrak{p}}}})/2 \to Pic(S_{0, \overline{k}}')/2\] 
is an isomorphism.  The Kummer sequence then shows that the right vertical arrow in (\ref{CD}) is an isomorphism.  What remains then is to find some extension $k'/k$ for which the base extension map $Br(S_{0}')[2] \to Br( S_{0, \overline{k}}')[2] \cong \zit/2$ is surjective (in this case,  $L/K_{\mathfrak{p}}$ will be the finite unramified Galois extension whose residue field extension is $k'/k$).  To this end,  we let $k'/k$ be a finite extension for which the action of $Gal(\overline{k}/k'),$ on $H^{2}_{\text{\'et}} (S_{0, \overline{k}}', \zit/2(1))$ is trivial.  Then,  note that the Hochschild-Serre spectral sequence:
\[ H^{p} (Gal(\overline{k}/k'),  H^{q}_{\text{\'et}} (S_{0, \overline{k}}', \zit/2(1))) \Rightarrow H^{p+q}_{\text{\'et}} (S_{0}', \zit/2(1))\]
degenerates,  since $k'$ is a finite field.  It follows that the base extension map 
\begin{equation} H^{2}_{\text{\'et}} (S_{0}', \zit/2(1)) \to H^{2}_{\text{\'et}} (S_{0, \overline{k}}', \zit/2(1))^{Gal(\overline{k}/k')} = H^{2}_{\text{\'et}} (S_{0, \overline{k}}', \zit/2(1))\label{last}\end{equation}
is surjective.  Since $S_{0}'$ is an Enriques surface,  $Br( S_{0, \overline{k}}') \cong \zit/2$.  In particular,  
\[ Br( S_{0, \overline{k}}')[2]^{Gal(\overline{k}/k')} = Br( S_{0, \overline{k}}') [2]\cong \zit/2 \]
So,  we deduce that the map $Br(S_{0}')[2] \to Br( S_{0, \overline{k}}')[2]$ is surjective,  as desired.
\end{proof}
\end{Lem}

\noindent Now,  let $\Gamma' :=  \pi_{\E'}^{*}A'\cup \pi_{\mathcal{S}'}^{*}B \in H^{3}_{\text{\'et}} (\W', \zit/2(2))$,  where $\E':= \E \times_{\mathfrak{o}} \mathfrak{o}'$ and $\W':= \W \times_{\mathfrak{o}} \mathfrak{o}'$.  This is a lift of $\gamma' \in H^{3}_{\text{\'et}} (W_{\overline{K_{\mathfrak{p}}}}, \zit/2(2))$.  We need to prove that $\gamma' \in H^{3}_{\text{\'et}} (W_{\overline{K_{\mathfrak{p}}}}, \zit/2(2))$ does not vanish under (\ref{func}).  To this end,  note that the specialization argument in \cite{G} shows that it is sufficient to prove that 
\[ \Gamma_{0}' \in H^{3}_{\text{\'et}} (W_{0,  \overline{k}}',  \zit/2(2))\]
the restriction of $\Gamma'$ to the geometric special fiber $W_{0,  \overline{k}}'$ of $\W' \to \text{Spec}(\mathfrak{o}')$ does not vanish under the map to Galois cohomology:
\begin{equation} H^{3}_{\text{\'et}} (W_{0,  \overline{k}}',  \zit/2(2)) \to H^{3} (F, \zit/2(2))\label{func2} \end{equation}
where $F=\overline{k}(W_{0,  \overline{k}}')$.  For this,  note that $W_{0}'$ is of the form $U \times S_{0}'$, where $U\cong \git_{m}$ (by Lemma \ref{Neron}) and $S_{0}'$ is the special fiber of $\mathcal{S}' \to \text{Spec}(\mathfrak{o}')$ as before.  Let $R \cong \pit^{1}_{k'}$ be the smooth completion of $U$ over $k'$ and consider $x \in R(\overline{k})$ not in $U(\overline{k})$.  Let $\nu$ be the valuation of $F$ corresponding to the divisor $x \times_{\overline{k}} S_{0,  \overline{k}}'$ on $R_{\overline{k}} \times_{\overline{k}} S_{0,  \overline{k}}'$ and let 
\[r_{\nu}:  H^{3} (F, \zit/2(2)) \to H^{2} (\overline{k}(S_{0,  \overline{k}}'), \zit/2(1))\]
be the corresponding residue.  Then,  the residue computation is the same as in \cite{CT} and \cite{G}.  Indeed,  we have
\begin{equation} r_{\nu} (\Gamma_{0}') = r_{x}(\alpha_{0}') \cdot \beta_{0} \in  H^{2} (\overline{k}(S_{0,  \overline{k}}'), \zit/2(1))\label{clinch}\end{equation}
where $\alpha_{0}' \in H^{1} (\overline{k} (R),  \zit/2(1))$ is the image of $A'$ under the composition 
\[ H^{1}_{\text{\'et}} (\E,  \zit/2(1)) \to H^{1}_{\text{\'et}} (E_{0,  \overline{k}},  \zit/2(1)) \to  H^{1} (\overline{k} (R),  \zit/2(1)),  \] 
$r_{x}: H^{1} (\overline{k} (R),  \zit/2(1)) \to \zit/2$ is the residue map corresponding to $x$ as above and $\beta_{0}$ is the image of $B$ under the composition 
\[ H^{2}_{\text{\'et}} (\mathcal{S}',  \zit/2(1)) \to H^{2}_{\text{\'et}} (S_{0,  \overline{k}}',  \zit/2(1)) \to  H^{2} (\overline{k}(S_{0,  \overline{k}}'),  \zit/2(1)) \]
The right hand side of (\ref{clinch}) is non-zero,  so long as $r_{x}(\alpha_{0}') \neq 0 \in \zit/2$.  There exist $x \in R(k')$ for which this holds; indeed,  we have
\[ \bigcap_{x \in R(\overline{k})} \text{ker } \{ r_{x}: H^{1} (\overline{k} (R),  \zit/2(1)) \to \zit/2\} = H^{1}_{\text{\'et}} (R_{\overline{k}}, \zit/2(1)) = 0 \]
Thus,  there is some $x$ for which (\ref{clinch}) is non-zero.  This shows that the image of $\Gamma_{0}'$ under (\ref{func2}) does not vanish,  as desired.

\section{An explicit example}
\noindent In this section, we give an example of $S$ as in the proof of Theorem \ref{main}.  There are conceivably many such examples, but we give one where all the assumptions may be checked directly. We need to check that the $S$ we construct satisfies Assumption \ref{ass}.  
\begin{Cons}\label{cons}
Let $X' \subset \pit^{3}$ be a quartic surface over $\qit$ for which there is a fixed-point-free automorphism of order $4$,  $\phi$,  defined over some number field.  We do not assume that $X'$ is smooth,  but we do assume the following:
\begin{enumerate}[label=(\alph*)]
\item\label{ext} $\phi$ extends to an action of $\pit^{3}$ with only finitely many fixed points,  and the corresponding fixed locus of $\iota:=\phi^2$ on $\pit^{3}$ consists of two skew lines;
\item\label{sing} the singular locus $X_{sing}'$ consists of a finite set of $\leq 8$ rational double points;
\item\label{inv} the fixed locus of $\iota$ on $X'$ contains $X_{sing}'$ (thus, the fixed locus is a scheme of length $8$ in which the smooth points of $X'$ have multiplicity $1$ and the double points muliplicity $2$).
\end{enumerate}
It is well-known that when \ref{sing} holds, the resolution $X$ obtained by blowing up the singular points is a $K3$ surface (this is also easy to check); $\phi$ then lifts to an action on $X$.  We also check in the lemma below that $\iota$ acting on $X$ via this induced action is a Nikulin involution.  In this case,  let $Y$ be the resolution of $X/\iota$ obtained by blowing up the $8$ rational double points of $X/\iota$; then, $Y$ is a K3 surface.  Moreover,  $\phi$ and $\phi^{3}$ both lift to fixed-point-free automorphisms on $X$,  which means that the induced action of $\phi$ on $Y$ is that of a fixed-point-free involution.  So, $S=Y/\phi$ is an Enriques surface.  Note that $S$ is defined over $\overline{\qit}$ (as is everything else in the construction).  Let $\pi: Y \to S$ be the corresponding double cover.  Morever,  let $\alpha \neq 0 \in H^{1} (S,  \zit/2) = Hom (\pi_{1} (S),  \zit/2)$ denote the corresponding class.
\end{Cons}
\begin{Lem} $\iota$ (acting on $X$) is a Nikulin involution.
\begin{proof} One may obtain $X$ as the strict transform of $X'$ relative to the blow-up $\epsilon: P\to \pit^{3}$ of $\pit^{3}$ at all the singular points of $X'$.  The action of $\phi$ on $\pit^{3}$ then lifts to $P$ and acts linearly on the exceptional divisors of $P$.  In particular,  the fixed locus of $\iota = \phi^{2}$ on each of these exceptional divisors is a union of projective linear subspaces,  and so $\iota$ does not act trivially on the exceptional divisors of $X' \to X$,  which are all conics.  The action of $\iota$ on $X$ thus fixes precisely $8$ points.  It now remains to show that any such involution is necessarily symplectic; or, equivalently,  that $Y$ is a $K3$ surface.  To this end,  we observe that $Y$ is certainly simply-connected.  So,  what remains is to show that the canonical divisor of $Y$ is trivial,  which is standard.  Indeed,  there is a branched double cover $g: \tilde{X} \to Y$,  where $\tilde{X}$ is the blow-up of $X$ along the $8$ fixed points of $\iota$.  Since $X$ is a $K3$ surface,  it follows that
\begin{equation} K_{\tilde{X}} = \sum_{i=1}^{8} E_{i} \in CH^{1} (\tilde{X})\label{exce} \end{equation}
where $E_{i}$ are the exceptional divisors of $\tilde{X}$.  We also have that 
\[ K_{\tilde{X}}= g^{*}K_{Y} + R \in CH^{1} (\tilde{X})\]
where $R$ is the ramification divisor of $g$.  However,  this latter is precisely the right-hand side of (\ref{exce}),  from which it follows that $g^{*}K_{Y} = 0 \in CH^{1} (\tilde{X})$.  Since $Y$ is simply-connected,  $CH^{1} (Y)$ is torsion-free, from which it follows that $K_{Y}  = 0 \in CH^{1} (Y)$
\end{proof}
\end{Lem}
\begin{Not}\label{not} We fix the following for the sequel:
\begin{itemize}
\item $p$: an odd prime;
\item $C'$: the quartic curve defined by $f(x_{0},  x_{1},  x_{2}):= x_{0}^{4} + px_{1}^{4} + x_{2}^{4} - 2x_{0}^{2}x_{2}^{2} -x_{0}x_{1}^{2}x_{2} =0$;
\item $X'$: the quartic surface defined by $g(x_{0},  x_{1},  x_{2},  x_{3}) = f(x_{0},  x_{1},  x_{2})+x_{3}^{4}=0$;
\item $\phi$: the automorphism on $\pit^{3}$ induced by $[x_{0}, x_{1}, x_{2}, x_{3}] \to [x_{0}, i x_{1}, - x_{2}, -i x_{3}]$,  where $i=\sqrt{-1}$.
\end{itemize}
\end{Not}
\begin{Lem}\label{prev} $X'$ and $\phi$ as in Notation \ref{not} satisfies the assumptions of Construction \ref{cons}. 
\begin{proof} One verifies directly that $X'$ is stable under the action of $\phi$ and that the set of fixed points of $\phi$ on $\pit^{3}$ is $\{ [1, 0, 0, 0],  [0,  1, 0,  0],  [0,  0,  1,  0],  [0,  0,  0,  1] \}$,  none of which lies on $X'$.  Additionally,  $\iota=\phi^{2}$ is given by 
\[[x_{0}, x_{1}, x_{2}, x_{3}] \to [x_{0}, -x_{1}, x_{2}, -x_{3}]\]
So,  the fixed locus of $\iota$ on $\pit^{3}$ is the  $L_{02} \cup L_{13}$ where $L_{ij}$ denotes the line in $\pit^{3}$ defined by $x_{i}=x_{j} =0$.  It follows that the fixed locus of $\iota$ on $X'$ is given by $F_{02} \cup F_{13}$, where 
\[ F_{ij} = \{ [x_{0},  x_{1},  x_{2},  x_{3}] \in \pit^{3} \ | \ x_{i} = x_{j} = g(x_{0},  x_{1},  x_{2},  x_{3}) = 0  \} \]
Note $\#F_{02} =4$ and $F_{13}$ consists of the two points $[1,  0,  \pm 1, 0]$.  As a finite scheme,  $F_{13} \subset \pit^{1}_{x_{0}x_{2}}$ is cut out by $(x_{0}-x_{2})^{2}(x_{0}+x_{2})^{2} = 0$.  Hence,  these two points each have multiplicity $2$.  Moreover,  a direct calculation shows that both of these points are singular points of $X'$,  which we note below are both rational double points.  What remains then is to establish that the only two singular points of $X'$ are $[1,  0,  \pm 1, 0]$.  Indeed,  a direct calculation shows that $X'$ is singular precisely at the points $[x_{0}, x_{1}, x_{2}, 0]$ where $[x_{0}, x_{1}, x_{2}]$ is a singular point of $C'$.  So,  we need to check that $[1,  0,  \pm 1]$ are the only singular points of $C'$.  Computing gradients of $f$,  the singular locus is defined by:
\begin{equation} 4x_{0}^{3} - 4x_{0}x_{2}^{2} -x_{1}^{2}x_{2} =4x_{2}^{3} - 4x_{0}^{2}x_{2} -x_{0}x_{1}^{2}= 4px_{1}^{3} -2x_{0}x_{1}x_{2}=0\label{singul} \end{equation}
If $x_{1} =0$,  we obtain the singular points $[1,  0,  \pm 1]$.  Moreover,  it is routine to check that these are rational double points on $C'$,  which implies that they are also rational double points on $X'$.  If $x_{1} \neq 0$,  the last equation in (\ref{singul}) gives $2px_{1}^{2} =x_{0}x_{2}$,  which leads to:
\[  4x_{0}^{3} - (4+\frac{1}{2p})x_{0}x_{2}^{2}=4x_{2}^{3} - (4+\frac{1}{2p})x_{0}^{2}x_{2}=0\]
One checks that the only solution to this system is $x_{0} =x_{2} = 0$,  which then forces $x_{1}=0$.  So,  $[1,  0,  \pm 1]$ are the only singular points.
\end{proof}
\end{Lem}

\noindent Since $X'$ now satisfies all the assumptions in Construction \ref{cons},  we let $Y$,  $S$,  $\gamma$,  $\pi$ be the associated objects corresponding to the $X'$ and $\phi$ in the construction. Additionally, note that $\phi$ acts on $C'$.  In fact, we have the following:
\begin{Lem}\label{prev2} The action of $\phi$ on the nodal curve $C'$ lifts to an action on its resolution $C$ in $X$ for which the induced action of $\iota=\phi^{2}$ is fixed-point free.  Moreover,  $C$ is a genus $1$ curve.
\begin{proof} Let $\rho: X \to X'$ be the blow-up of $X'$ along the two singular points of $X'$.      This resolves the singularities of $X'$, as was observed in the proof of Lemma \ref{prev}.  Additionally, $C = \rho^{-1} (C)$ is a smooth curve; a computation shows that its genus is $1$. The action of $\phi$ on $C'$ lifts to an automorphism on $C$,  since the latter is obtained by blowing up $C$ on what happen to be the set of fixed points of $\iota$ (acting on $C'$).  What remains is to see that the induced action of $\iota=\phi^{2}$ on $C$ is fixed-point-free.  For this,  note that since $\iota$ acts on $C$ as an involution,  either $\iota$ is fixed-point-free on $C$ or the quotient $C/\iota$ has genus $0$.  However, the latter case cannot occur,  since $C/\iota$ is also stable under the action of the fixed-point-free involution $\phi$ acting on $Y$ (impossible if $C/\iota$ has genus $0$).  It follows that the action of $\iota$ on $C$ is fixed-point-free, as desired.
\end{proof}
\end{Lem}
\noindent As observed in the proof of Lemma \ref{prev2},  $\tilde{E}:=C/\iota \subset Y$ is a genus $1$ curve.  Since $\pi$ is \'etale, $E = \pi(\tilde{E})$ is also a genus $1$ curve on $S$.  Also observe that $\tilde{E},  E$ are defined over $K=\qit(i)$.  Note that since $C(K) \neq \emptyset$,  the same is true of $\tilde{E}$. Thus, let $\tilde{E},  E/K_{p}$ denote the above genus $1$ curves, viewed as curves over $K_{\mathfrak{p}}$,  the completion of $K$ at $\mathfrak{p}$,  a prime lying over $p$.  What remains now is to verify that this construction satisfies the reduction hypothesis in Assumption \ref{ass}:
\begin{Lem} $\tilde{E}$ has multiplicative reduction at a prime $\mathfrak{p}$ lying over $p$.
\begin{proof}
Let $\C'$ denote the subscheme of $\pit^{2}_{\mathfrak{o}}$ defined by 
\[ f(x_{0},  x_{1},  x_{2}):= x_{0}^{4} + px_{1}^{4} + x_{2}^{4} - 2x_{0}^{2}x_{2}^{2} -x_{0}x_{1}^{2}x_{2} =0 \] 
The generic fiber of $\C' \to \text{Spec} (\mathfrak{o})$ is certainly $C'$ and the special fiber $C_{0}'$ is the curve in $\pit^{2}_{k}$ (where $k$ denotes the residue field) defined by 
\[ (x_{0}^{2} - x_{2}^{2})^{2} -x_{0}x_{1}^{2}x_{2} =0\]
A direct gradient computation shows that the singular locus of $C_{0}'$ consists of $[1,  0,  \pm 1]$ and $[0, 0, 1]$; these also happen to be exactly the fixed points of $\iota$ on $C_{0}'$.  
The points $[1,  0,  \pm 1] \in \C'(K_{\mathfrak{p}})$ may be viewed as sections $s_{\pm}:\text{Spec} (\mathfrak{o}) \to \C'$.  By abuse of notation,  let $s_{\pm}$ also denote the corresponding closed subschemes of $\C'$.  Blowing up along $s_{+} \cup s_{-}$ yields a scheme $\C \to \text{Spec} (\mathfrak{o})$,  whose generic fiber is certainly $C$.  The involution $\iota$ on $C$ extends to an action on $\C \to \text{Spec} (\mathfrak{o})$.  The special fiber $C_{0}$ is a nodal rational curve (since $C_{0}'$ had $3$ nodes).  Let $C_{0,  sm} \subset C_{0}$ be the complement of this node; i.e.,  the smooth part of $C_{0}$.  Then,  $C_{0,  sm} \cong \git_{m}$  By construction,  the action of $\iota$ on $C_{0,  sm}$ is fixed-point-free.  It follows that the quotient $C_{0,  sm}/\iota$ is also isomorphic to $\git_{m}$,  which gives the desired multiplicative reduction.

\end{proof}
\end{Lem}

\section{A remark on the proof of Theorem \ref{main}}

\noindent It is well-known that if an Enriques surface $S$ contains genus $1$ curves,  then there is an elliptic pencil $h: S \to \pit^{1}$.  In particular,  one may ask whether we could have instead taken the approach of \cite{BO} and \cite{CT},  using the fact that $S$ (in the proof of Theorem \ref{main}) admits an elliptic fibration, $h: S \to \pit^{1}$.  This argument does not work out,  as it turns out.  Indeed,  let $\alpha$ be the class from the proof of Theorem \ref{main},  corresponding to the double cover $Y \to S$ and let $E_{t} =h^{-1}(t)$.  The claim is that $\restr{\alpha}{E_{t}} \in H^{1}(E_{t} ,  \zit/2)$ vanishes for almost all $t \in \pit^{1}(\cit)$.  This will follow if we show that the Enriques involution $\phi$ that acts on $Y$ acts non-trivially on the pencil $|\tilde{E}|$ (in the notation of the proof of Theorem \ref{main}),  so that the $\pi^{-1}(E_{t})$ is a disjoint union of two fibers in the pencil $|\tilde{E}|$.  To this end,  note that $\phi$ acts on $\tilde{E} \subset Y$ by assumption and,  hence,  also acts on the pencil $|\tilde{E}|$.  If $\phi$ were to act trivially on $|\tilde{E}|$,  then it would act on the corresponding elliptic pencil $\tilde{h}: Y \to \pit^{1}$.  Let $U \subset \pit^{1}$ be the open subscheme over which $\tilde{h}$ is smooth.  Then,  $\phi$ would act on $Y_{U}:= \tilde{h}^{-1}(U) \to U$.  Moreover,  since $\phi$ acts freely,  it would then follow that it acts by translation on the fibers of this latter.  In particular,  $\phi^{*}$ would act trivially on $H^{1}(U,  R^{1}\tilde{h}_{*}\cit(1))\subset  H^{2} (Y_{U},  \cit(1))$.  Now,  the Gysin sequence gives an injective map:
\[ \cit\cong H^{2,0} (Y) \to H^{1}(U,  R^{1}\tilde{h}_{*}\cit(1)) \]
via restriction.  This would imply that $\phi^{*}$ acts trivially on $H^{2,0} (Y)$,  which would mean $\phi^{*}$ acted symplectically.  This gives the desired contradiction,  from which it follows that $\phi$ acts non-trivially on $|\tilde{E}|$.  Thus,  $\phi$ fixes exactly two of the curves in $|\tilde{E}|$,  one of which is $\tilde{E}$ itself.  One can show that the corresponding pencil on $S$ then has a multiple fiber whose corresponding reduced curve is $E$.  Since $\phi$ acts non-trivially on $|\tilde{E}|$,  one can likely show that the $Y$ and $S$ considered here are cases of the more general examples considered in \cite{HS}.

\end{document}